\documentclass[10pt]{article}
\usepackage{amsmath,amsfonts,amsthm,amssymb,amscd}

\catcode`\é=\active\def é{\'e} \catcode`\è=\active\def è{\`e}
\catcode`\ç=\active\def ç{\c{c}} \catcode`\à=\active\def à{\`a}
\catcode`\ê=\active\def ê{\^e} \catcode`\ù=\active\def ù{\`u}
\catcode`\ô=\active\def ô{\^o} \catcode`\û=\active\def û{\^u}
\catcode`\î=\active\def î{\^{\i}} \catcode`\â=\active\def â{\^a}
\catcode`\ö=\active\def ö{\"o}

\newcommand{\be}{\begin{equation}}
\newcommand{\ee}{\end{equation}}
\newcommand{\bs}{\begin{split}}
\newcommand{\es}{\end{split}}
\newcommand{\ba}{\begin{align}}
\newcommand{\ea}{\end{align}}
\newcommand{\basl}[1]{\begin{align}\begin{split}\label{#1}}
\newcommand{\bas}{\begin{align}\begin{split}}

\newtheorem{theo}{Theorem}[section]

\newtheorem{coro}[theo]{Corollary}
\newtheorem{rem}[theo]{Remark}

\newcommand\fpr{\hfill$\Box$\null}

\newcommand\N{\mathbb{N}}

\newcommand\R{\mathbb{R}}
\newcommand\C{\mathbb{C}}

\newcommand\ce{{\cal E}}
\newcommand\cf{{\cal F}}

\newcommand\cl{{\cal L}}

\newcommand\ga{\gamma}

\newcommand\si{\sigma}

\newcommand\la{\lambda}

\newcommand\al{\alpha}

\def\fr #1#2{\frac{#1}{#2}}

\title{Inverse spectral results for Schr\"odinger operators on the unit interval with potentials in $L^P$ spaces}
\author{L. Amour\footnote{laurent.amour@univ-reims.fr}\ \ and T. Raoux\footnote{thierry.raoux@univ-reims.fr}}
\date{September 27, 2007}

\begin{document}
\maketitle \centerline{Laboratoire de Math\'ematiques EDPPM,
UMR-CNRS 6056} \centerline{ Universit\'e de Reims,
 Moulin de la Housse}
 \centerline{BP 1039,
 51687 REIMS Cedex 2, France.}
\vskip 1cm
\begin{abstract}
We consider the Schr\"odinger operator on $[0,1]$ with potential in
$L^1$. We prove that two potentials already known on $[a,1]$
($a\in\left(0,\frac{1}{2}\right]$) and having their difference in
$L^p$ are equal if the number of their common eigenvalues is
sufficiently large. The result here is to write down explicitly this
number in terms of $p$ (and $a$) showing the role of $p$.
\end{abstract}

% \tableofcontents

\parindent=0pt
\parskip 10pt
\baselineskip 15pt

\section{Introduction}

This paper is concerned with the Schr\"odinger operator, \be
A_{q,h,H}=-\fr{d^2}{dx^2}+q \ee on $[0,1]$ associated with the
boundary conditions, \be u'(0)+hu(0)=0,\ \quad u'(1)+Hu(1)=0.\ee
Here the potential $q$ is a real-valued function belonging to
$L^1([0,1])$ and $h,H\in \R$. The corresponding spectrum is a
sequence of simple eigenvalues. Let us denote by
$(\la_j(q,h,H))_{j\in\N\cup\{0\}}$ the increasing sequence of its
eigenvalues. Then, the asymptotic expansion is ([LG]) \be\label{ae}
\la_j(q,h,H)=j^2\pi^2\,+\,2(H-h)\,+\,\int_0^1\,q(x) dx\,+\,
o(1)\quad{\rm as}\ j\rightarrow +\infty.\ee

In 1978, H. Hochstadt and B. Lieberman [HL] proved that the whole
spectra of $A_{q,h,H}$ determine uniquely $q$ when it is already
known on $[\frac{1}{2},1]$. More precisely, if $q_1,q_2\in
L^1[0,1]$, $q_1=q_2$ a.e. on $[\frac{1}{2},1]$ and if the spectrum
of $A_{q_1,h_1,H}$ is exactly the spectrum of $A_{q_2,h_2,H}$ then
$q_1=q_2$ and $h_1=h_2$. In 2000, several extensions of this result
are given by F. Gesztesy and B. Simon [GS] considering that the
potentials belonging to $L^1([0,1])$ are known on a larger interval
($[a,1]$ with $a\in (0,\frac{1}{2}]$) and assuming that the common
spectrum of $A_{q_1,h_1,H}$ and $A_{q_2,h_2,H}$ is sufficiently
large (in term of $a$). Another result in [GS] is to assume that the
potential belongs to $C^k$ so that the number of common eigenvalues
is given in term of $(a,k)$. Our aim here is to obtain a similar
result for potentials in $L^p$. Actually only their difference needs
to be in $L^p$. We shall derive that two potentials already known on
$[a,1]$ and having their difference in $L^p$ are equal if their
common spectrum is sufficiently large (depending of $(a,p)$) (see
theorem \ref{t1} below). Let us also mention at this point that our
proof is different from the proof in [GS]).

For any $\al=(\al_j)_{j\in\N\cup\{0\}}$, $\al_j\in\C$, set
$$n_\al(t)=\sharp\{j\in\N\cup\{0\}\quad |\quad \vert \al_j\vert\leq t\},\ \forall\,t\geq
0.$$

Our purpose here is to prove the following result
\begin{theo}\label{t1}
Fix $q_1,q_2\in L^1([0,1])$ and $h_1,h_2,H\in \R$. Consider the
infinite set $S$ \be S\subset \si(A_{q_1,h_1,H})\cap
\si(A_{q_2,h_2,H})\ee

Fix $a\in (0,\fr{1}{2}]$ and $p\in [1,+\infty)$. Suppose that
$q_1=q_2$ on $[a,1]$ and $q_1-q_2\in L^p([0,a])$. Assume that there
exists a real number $C$ such that \[ 2a\,n_{\si(A)}(t)+C\ \geq\
n_S(t)\ \geq\ 2a\, n_{\si(A)}(t)+\fr{1}{2p}-2a,\ t\in S,\ t\ large\
enough,\quad(H1)
\] where $A$ denotes either $A_{q_1,h_1,H}$ or $A_{q_2,h_2,H}$. Then
$h_1=h_2$ and $q_1=q_2$.
\end{theo}
Roughly speaking, theorem \ref{t1} says that the potential given on
$[a,1]$ together with a sufficiently large part (depending on $a$)
of its spectrum determine entirely the potential on $[0,1]$.

In the particular case $p=1$, a similar result (among many others)
is proved in [GS, theorem 1.3] with the following modifications: in
[GS],

$(i)$ the lower bound in $(H_1)$ is $2a\, n_{\si(A)}(t)+\fr{1}{2}-a$
and $t\in\R$, $t>0$ is large enough.

$(ii)$ there is no upper bound in $(H_1)$.

Concerning $(ii)$ our result is weaker and concerning $(i)$ it is
stronger because of the following two reasons: 1. we have $-2a$
instead of $-a$. 2. the parameter $t$ needs only to be in $S$. These
two points are involved in proof of the next corollary. \newline The
upper bound in $(H_1)$ imposes here that the given spectrum is in
some sense regularly spaced. This is not required in [GS]. However,
because of the above points 1 and 2 our lower bound is well-adapted
to results like : the even spectrum and the potential given on
$\left[\frac{1}{4},1\right]$ determine the potential on $[0,1]$ (see
corollary below) whereas in that case, [GS] needs slightly more than
half of the spectrum (see remark below theorem 1.3 in [GS]).

\begin{coro}
The even (resp. odd) spectrum $(\la_{2j}(q,h,H))_{j\geq 0}$ (resp.
$(\la_{2j+1}(q,h,H))_{j\geq 0}$) and $q|_{[0,\fr{1}{4}]}$ determine
$q$ on $[0,1]$
\end{coro}

{\it Proof :} For the even (resp. odd) case, apply theorem \ref{t1}
with $a=\fr{1}{4}$, $S=(\la_{2j})_{j\geq 0}$ (resp.
$S=(\la_{2j+1})_{j\geq 0}$) and use $n_{\si(A)}(\la_{2j})=2j+1$ and
$n_{S}(\la_{2j}))=j+1$ (resp. $n_{\si(A)}(\la_{2j+1}))=2j+2$ and
$n_{S}(\la_{2j+1}))=j+1$).\fpr

\begin{rem}
Similar results may be obtained also for the Dirichlet boundary
conditions. Moreover, this method may be applied analogously to the
AKNS systems. For AKNS systems one may also refer to the work in
[DG] where the given spectrum is regularly spaced : $S=\{\la_{jk},\
j\geq 0\}$ with $k$ being a fixed positive integer.
\end{rem}

Whereas the proof of the results in [GS] relies on the Weyl-Titmarsh
functions, the starting point here is different and it is based on
an idea taken in ([L]) (which appears in a short proof that two
spectra determine the potential). Let us describe the main points.
1. An entire function $f$ depending on $q_1-q_2$ restricted to
$[0,a]$ is introduced having the property to vanish on the common
eigenvalues (to be complete, we mention that $f(z)$ is the r.h.s. of
(\ref{ep})). 2. The second step is to use the growth property of $f$
to derive that it is identically vanishing which directly follow
from the maximum modulus principle. 3. The last step is to derive
that $f\equiv 0$ implies that $q_1\equiv q_2$.

Our contribution here is to modify the second step above in order to
deal with potentials in $L^p$. The main fact is to replace the
maximum modulus principle by a result of Levinson stated in Levin
[L] (see Step 4 below). Whereas the maximum modulus principle is
applied to $f$, we shall apply Levinson's result to the Fourier
transform of $f$. More precisely, we rather use the Fourier
transform of $f$ that we call $g$ for the following two reasons. The
first one is that $f$ is actually roughly speaking close to the
inverse Fourier transform of $q_1-q_2$, so that $g$ is close in some
sense to $q_1-q_2$ and it is expected that this imply that the
assumption $(q_1-q_2)\in L^p$ is rewritten as $g\in L^p$ without any
loss of information. The second one is (since the inverse Fourier
transform of $g$ vanishes on the common eigenvalues) to remark a
result given in [L] due to Levinson and essentially stating that the
inverse Fourier transform of function being in $L^p$ is entirely
vanishing if it has a sufficiently large number (depending on $p$)
of zeros. Therefore, this shall replace in our proof the point 2
above.

This work is concerned with $L^p$ spaces and a work involving others
spaces is in progress.

In the next section we establish theorem \ref{t1}. Its proof is
split into 5 steps. The first step is to define properly the
function $g$ and to give some of the properties that shall be used
in the sequel. The second point is to recall that the inverse
Fourier transform of $g$ is vanishing on the common spectra of
$A_{q_1,h_1,H}$ and $A_{q_2,h_2,H}$. In the third step we introduce
an auxiliary property $(H_2)$ derived from $(H_1)$. The fourth step
consists in proving that if $(H_2)$ is satisfied then $g$ is
vanishing. It is at this point that we use Levinson's result. In the
fifth step we give short proof of the already known fact: $g$ equals
zero implies $q_1=q_2$ and $h_1=h_2$.

\section{Proof of theorem \ref{t1}}

\subsection{Step 1: definition of $g$}

{\it $\bullet$ Definition of $\psi$:} For $z\in\C$, let
$\psi(\cdot,z,q,h)$ defined on $[0,1]$, be the solution to
$(-\fr{d^2}{dx^2}+q)\psi=z\psi$, $\psi(0)=1$, $\psi'(0)=-h$. It is
known that $\psi(x,\cdot,q,h)$ is an entire function and ([LG]) \be
\psi(x,z,q,h)=\cos\sqrt{z}x + O\left(\fr{e^{|\Im
\sqrt{z}|x}}{\sqrt{|z|}}\right),\quad {\rm as}\ |z|\rightarrow
+\infty,\ee uniformly in $x\in [0,1]$.

{\it $\bullet$ Definition of $r$:} Fix $h_1,h_2\in\R$ and
$q_1,q_2\in L^1([0,1])$. For $x\in [0,1]$, $z\in\C$, let
$r(z,x)=-\psi(x,z^2,q_1,h_1)\psi(x,z^2,q_2,h_2)+\fr{1}{2}(1+\cos
2{z}x)$. Clearly, \be r(z,x)= O\left(\fr{e^{2|\Im z|x}}{|z|}\right),
\quad {\rm as}\ |z|\rightarrow +\infty,\ee uniformly in $x\in
[0,1]$.

In order to apply Levinson's result on $[-2a,2a]$ in step 4 we
introduce below a scaling and extension by parity operator $\cal E$
which shall be always applied to $r$ and $q_1-q_2$ in the sequel. In
particular it allows to define below $s$ with the usual Fourier
transform instead of the cosine Fourier transform.

 {\it
$\bullet$ Definition of $\ce$:} For any real-valued function $u$
defined a.e. on $[0,a]$, let us set
$\ce(u)(x)=u\left(\fr{|x|}{2}\right)$, for a.e $x\in[-2a,2a]$. If
$v$ is real-valued function depending on two variables with the
second variable belonging to $[0,a]$, then $\ce_2(v)$ denotes $\ce$
applied to the second variable. In particular, \be\label{1}
\ce_2(r)(z,x)= O\left(\fr{e^{|\Im z||x|}}{|z|}\right),\quad {\rm
as}\ |z|\rightarrow +\infty,\ee uniformly in $x\in[-2a,2a]$.

{\it $\bullet$ Definition of $s$:} Set
$s(\cdot,x)=\cf_1\ce_2(r)(\cdot,x)$ on $\R$,
$\forall\,x\in[-2a,2a]$. Here $\cf_1$ is the Fourier transform
applied on the first variable. We remark that (Paley-Wiener
theorem), \be\label{2} {\rm supp}\,s=\{(y,x)\in[-2a,2a]^2\ |\
|y|\leq|x|\} \ee and using (\ref{1}) and the regularity properties
of $r$,\be\label{3} \forall\, q>1,\ \sup_{x\in[-2a,2a]}\Vert
s(\cdot,x)\Vert_{L^q[-2a,2a]} <\infty.\ee

{\it $\bullet$ Definition of $g$:} Applying successively Hölder
inequality, Fubini's theorem and Cauchy-Schwarz inequality show
$$
\int_{-2a}^{2a}\left|\int_{-2a}^{2a}s(y,x)\ce(q_1-q_2)(x)dx\right|^p
dy$$
\begin{eqnarray}\label{HFCS}
&\leq& \left(\int_{-2a}^{2a}| \ce(q_1-q_2)(x)|\,dx\right)^{p-1}
\int_{-2a}^{2a}\int_{-2a}^{2a}\left|s(y,x)\right|^p\left|\ce(q_1-q_2)(x)\right|\,dxdy\nonumber\\
\\
&\leq& (2^{p+1}\sqrt{a}) \Vert q_1-q_2\Vert_{L^1([0,a])}^p
\sup_{x\in[-2a,2a]}\Vert
s(\cdot,x)\Vert_{L^{2p}[-2a,2a]}^p.\nonumber\end{eqnarray} Then
$(\ref{HFCS})$ yields that $g$ is well-defined by, \be\label{4}
g(y)= \ce(q_1-q_2)(y)-2\int_{-2a}^{2a} s(y,x)\ce(q_1-q_2)(x)\,dx,
{\rm \ for\ a.e.}\ y\in [-2a,2a].\ee Since $q_1-q_2\in L^p([0,a])$
then \be \label{glp} g\in L^p([-2a,2a]). \ee We shall also use the
notation $T_s$ for the integral operator with kernel $2s$. In
particular,
$$T_s:\ L^p([-2a,2a])\rightarrow L^p([-2a,2a])$$ and \be\label{5}
g=(1-T_s)(\ce(q_1-q_2)). \ee
\begin{rem}\ \newline
It may be also natural to consider for $z\in\C$, $f(z)=\int_0^a
\left(-1
 +2\psi(x,z^2,q_1,h_1)\psi(x,z^2,q_2,h_2)\right)(q_1(x)-q_2(x))dx$
 (see (\ref{ep}))
 and to define $g$ as the Fourier transform of $f$. However it is
 not directly clear that $g$ would be in $L^p$ when $q_1-q_2$ is
 $L^p$, in particular $g(y)$ may be not written as $\int_{\R}
 e^{-izy} f(z)dz$.  Moreover $f$ itself is not involved in the main step (step 4).
 Therefore we choose to define $g$ by (\ref{4}),
 verify (\ref{glp})(\ref{5}) and check that ${\cal F}^{-1}g=f$ (cf
 (\ref{ep})).
\end{rem}

\subsection{Step 2: $z^2\in S\Rightarrow\int_{-2a}^{2a}
e^{izy}g(y)\,dy=0$} Since $g\in L^1([-2a,2a])$ then (\ref{4})
implies that $\int_{-2a}^{2a}
e^{izy}g(y)\,dy=\cf^{-1}(\ce(q_1-q_2))(z)-2\int_{-2a}^{2a}\int_{-2a}^{2a}e^{izy}s(y,x)\ce(q_1-q_2)(x)dxdy.$
\newline $\int_{-2a}^{2a} e^{izy}g(y)\,dy$
\begin{eqnarray}\label{ep}
 &=&\cf^{-1}(\ce(q_1-q_2))(z)-2\int_{-2a}^{2a}\int_{-2a}^{2a}e^{izy}\cf_1\ce_2(r)(y,x)\ce(q_1-q_2)(x)dxdy\nonumber\\
 &=&\int_{-2a}^{2a}(e^{izx}-2\ce_2(r)(z,x))\ce(q_1-q_2)(x)dx\nonumber\\
 &=&4\int_{0}^{a}(\cos 2zx-2r(z,x))(q_1-q_2)(x)dx\nonumber\\
 &=&4\int_0^a \left(-1
 +2\psi(x,z^2,q_1,h_1)\psi(x,z^2,q_2,h_2)\right)(q_1(x)-q_2(x))dx.
\end{eqnarray}
This is known to vanish for $z\in S$. Indeed, multiply
$\left(-\fr{d^2}{dx^2}+q_1(x)-z^2\right)\psi(x,z^2,q_1,h_1)=0$ by
$\psi(x,z^2,q_2,h_2)$, multiply
$\left(-\fr{d^2}{dx^2}+q_2(x)-z^2\right)\psi(x,z^2,q_2,h_2)=0$ by
$\psi(x,z^2,q_1,h_1)$ and integrate their difference on $[0,1]$ to
obtain that the term in (\ref{ep}) equals $2(h_1-h_2)+\int_0^1
q_2(x)-q_1(x) dx$. This term is zero from $(\ref{ae})$ since it is
assumed that $S$ contains an infinite number of points.

\begin{rem}\label{r2}
Since the translation on the potential $q$ acts as a translation on
the spectrum of $A_{q,h,H}$, it is assumed without loss of
generality that the $\la_j(q_1,h_1,H)$'s and the
$\la_j(q_2,h_2,H)$'s are positive real numbers.
\end{rem}

Let us write $S=(s_j)_{j\geq 0}$ where $(s_j)$ is an increasing
sequence and consider the following property: there exists $D\in\R$
such that,
\[ \quad\fr{\pi}{2a}j+D\ \leq\ \sqrt{s_j}
\ \leq\
\fr{\pi}{2a}\left(j+1-\fr{1}{2p}\right)+O\left(\fr{1}{j}\right),\quad
{\rm as} \ j\rightarrow\infty\quad\quad (H_2)\]

\subsection{Step 3: $(H_1)\Rightarrow (H_2)$}
Since $s_j\in(\la_n(q_1,h_1,H))_{n\in\N\cup\{0\}}$ then
$s_j=\la_{m_j}(q_1,h_1,H)$ for some increasing sequence of integer
numbers $(m_j)_{n\in\N\cup\{0\}}$. In particular,
$n_{\si(A_{q_1,h_1,H})}(s_j)=m_j+1$ and following $(\ref{ae})$ we
have
\begin{eqnarray}\label{9}
 n_{\sigma(A_{q_1,h_1,H})}(s_j)&=&\sharp\{n\geq 0\ |\ \la_n(q_1,h_1,H)\leq s_j\}\nonumber\\
 &=&\sharp\left\{n\geq 0\ |\ n\pi+O\left(\fr{1}{n}\right)\leq
 m_j\pi+O\left(\fr{1}{m_j}\right)\right\}.
\end{eqnarray}
We also have \be\label{10}n_S(s_j)=j+1.\ee Let us check that
$(H_1)$(\ref{9})(\ref{10}) imply $(H_2)$.

On one side, the second inequality in $(H_1)$ together with
(\ref{9})(\ref{10}) read as
\[j+1\geq
2a\left(\fr{\sqrt{s_j}}{\pi}+1+O\left(\fr{1}{j}\right)\right)+\fr{1}{2p}-2a\]
which is the second inequality in $(H_2)$. On the other side, the
first inequality in $(H_1)$ together with (\ref{9})(\ref{10}) give
$j+1\leq
2a\left(\fr{\sqrt{s_j}}{\pi}+1+O\left(\fr{1}{j}\right)\right)+C$,
that is to say,
\[
\sqrt{s_j}\geq
\pi\left(\fr{j+1-C}{2a}-1\right)+O\left(\fr{1}{j}\right),
\]
which imply the first inequality in $(H_2)$.

\subsection{Step 4 : $(H_2)\Rightarrow g=0$}

Let $\al=(\al_j)$ be a sequence of complex numbers and define
$N_\al(R)=\int_0^R\fr{n_\al(t)}{t}\,dt$. Let us recall the following
result [L, app. III, sec. 2, th. 3]:

{\it Let $q\in]1,+\infty]$, $\al=(\al_j)$ a sequence of complex
numbers satisfying $$ \limsup_{R\rightarrow\infty}
N_\al(R)-2R+\fr{1}{q}\ln R >-\infty.$$
 Then the family $\{y\mapsto
e^{i\al_j y}\}$ is complete in $L^q([-\pi,\pi])$ (if $q=\infty$ then
$\{y\mapsto e^{i\al_j y}\}$ is complete in $C^0([-\pi,\pi])$.}

Set (see remark \ref{r2})
\[S^\fr{1}{2}=\{\pm\sqrt{s_j},\ j\geq 0\},\quad S^{\fr{1}{2},+}=\{\sqrt{s_j},\ j\geq
0\}\]

In particular, replacing $\pi$ by $2a$ and setting $\al=S^\fr{1}{2}$
this result read as: if$$ \limsup_{R\rightarrow\infty}
N_{S^\fr{1}{2}}(R)-\fr{4a}{\pi}R+\fr{1}{q}\ln R
>-\infty\qquad (H_L)$$ then $\{e^{i\al_j y}\}$ is complete in $L^q([-2a,2a])$
(resp. $C^0([-2a,2a])$) if $q\in ]1,+\infty[$ (resp. if
$q=+\infty$). To see this one may either reproduce the proof taken
in [L] while replacing $\pi$ with $2a$ or use Riesz representation
theorem to write $L\in (L^q([-2a,2a]))'$ as
$L(h)=\int_{-2a}^{2a}l(y)h(y)dy$ (for some $l\in L^{q'}([-2a,2a])$)
and make the change of variable $y'=\frac{\pi}{2a}y$, then use
$n_{\al}(t)=n_{\la\al}(\la t)$ and
 $N_{\al}(t)=N_{\la\al}(\la t)$ for all $\la>0$ and $t>0$.

Let us check that \be\label{20} (H_2)\quad \Longrightarrow
(H_L)\qquad{\rm if}\quad \frac{1}{p}+\frac{1}{q}=1.\ee Since $s_j>0$
for all $j\geq 0$, $n_{S^\fr{1}{2}}(t)=2n_{S^{\fr{1}{2},+}}(t)$ for
all $t>0$. Moreover, $n_{S^{\fr{1}{2},+}}(t)=j+1,\quad\forall\,t\in
[\sqrt{s_j},\sqrt{s_{j+1}}),\ \forall\, j\geq 0$. Therefore,
$N_{S^{\fr{1}{2}}}(\sqrt{s_j})=2\sum_{k=0}^{j-1}
\int_{\sqrt{s_k}}^{\sqrt{s_{k+1}}} \fr{k+1}{t}\,dt$ which gives
\be\label{21} N_{S^{\fr{1}{2}}}(\sqrt{s_j})=2\left(j\ln
\sqrt{s_j}-\sum_{k=0}^{j-1}\ln \sqrt{s_k}\right),\quad\forall\,
j\geq 0\ee

Following $(H_2)$ (left inequality), \be\label{22} j\ln
\sqrt{s_j}\geq j\left(\ln j+\ln\fr{\pi}{2a}\right)+O(1),\quad{\rm
as}\ j\rightarrow +\infty.\ee Following $(H_2)$ (right inequality)
and $\ln j!=\left(j+\fr{1}{2}\right)\ln j -j+O(1)$ as $j\rightarrow
+\infty$,
\begin{eqnarray}\label{23}
\sum_{k=0}^{j-1}\ln
\sqrt{s_k}&\leq&j\ln\fr{\pi}{2a}+\sum_{k=0}^{j-1} \ln (k+1) +
\ln\left(1-\fr{\fr{1}{2p}+O\left(\fr{1}{k}\right)}{k+1}\right)\nonumber\\
&\leq&j\ln\fr{\pi}{2a}+ \ln
j!-\sum_{k=0}^{j-1}\fr{\fr{1}{2p}+O\left(\fr{1}{k}\right)}{k+1}+
O\left(\fr{1}{k^2}\right)\\
 &\leq&j\ln\fr{\pi}{2a}+\left(j+\fr{1}{2}\right)\ln j -j-\fr{1}{2p}(\ln j
 +\ga)+O(1)\nonumber
\end{eqnarray}
where $\gamma$ is the Euler constant and as $j\rightarrow +\infty$.
In particular, combining (\ref{21}) with (\ref{22})(\ref{23}) and
using again $(H_2)$ give
\begin{eqnarray}
N_{S^\fr{1}{2}}(\sqrt{s_j})&\geq& \left(\fr{1}{p}-1\right)\ln
j +2j+O(1)\nonumber\\
 &\geq& \left(\fr{1}{p}-1\right)\ln
 \sqrt{s_j}+\fr{4a}{\pi}\sqrt{s_j}+O(1)
\nonumber
\end{eqnarray}
as $j\rightarrow\infty$. This proves (\ref{20}) and $\{e^{\pm
i\sqrt{s_j} y}\}$ is complete in $L^q([-2a,2a])$ (resp.
$C^0([-2a,2a])$) if $q\in ]1,+\infty[$ (resp. if $q=+\infty$). For
$p\in ]1,+\infty[$ (resp. $p=1$) define $L\in (L^q([-2a,2a]))'$
(resp. $L\in (C^0([-2a,2a]))'$ by
$$
L:\left\{\begin{array}{ccl} L^q([-2a,2a])\ ({\rm resp.}\
C^0([-2a,2a]) &\rightarrow&\C\\
h&\mapsto& \int_{-2a}^{2a} h(y)g(y)\,dy
\end{array}\right..
$$
According to step 1, $L$ is well-defined. Following step 2,
$L(y\mapsto e^{\pm i\sqrt{s_j} y})=0$ for all $j$. The completeness
property of $\{y\mapsto e^{\pm i\sqrt{s_j} y}\}$ implies that $L$ is
vanishing identically. This proves that $g\equiv 0$.

\subsection{Step 5 : $g=0\Rightarrow (h_1,q_1)=(h_2,q_2)$}

The fact that $h_1=h_2$ shall follow the asymptotic expansions of
the eigenvalues and the fact that $q_1=q_2$ is actually proved in
[L] since $g=0$ implies $\int_0^1 (-\fr{1}{2}
 +\psi(x,z,q_1,h_1)\psi(x,z,q_2,h_2)(q_1-q_2)(x)dx=0,\ \forall
 z\in\R$. For sake of completeness, let us give a
 shorter proof involving only the function $s$ and its properties. Indeed,
 (\ref{2})(\ref{3}) and the definition of $T_s$,
 Fubini's theorem yield for any $u\in L_{\rm }^1([-2a,2a])$
  and for any $n\geq 1$
 $$\Vert
T_s^n u\Vert_{L_{\rm }^1([-2a,2a])}$$ $$ \leq
2^n\int_{\left\{|t_{n+1}|\leq\dots\leq |t_{2}|\leq|t_{1}|\leq
2a\right\}}|s(t_{n+1},t_n)\cdots s(t_2,t_1)u(t_1)|\,dt_1\dots
dt_{n+1}.$$

Consequently $$ \Vert T_s^n \Vert_{\cl(L^1([-2a,2a]))}\leq \fr{(4
\sqrt a \sup_{x\in[-2a,2a]}\Vert
s(\cdot,x)\Vert_{L^2[-2a,2a]})^n}{\sqrt{n!}}.
$$ In particular, this proves that $T_s^n$ is a contracting map in $L^1([-2a,2a])$ for $n$ large enough.
It follows that $0$ is the unique fixed point of $T_s$ in
$L^1([-2a,2a])$. By (\ref{5}), this proves that $q_1=q_2$ and the
fact that $h_1=h_2$ follows the asymptotic expansion of the
eigenvalues (\ref{ae}).

\vskip1cm \centerline{\bf References}

[DG] R. del Rio and B. Gr\'ebert, {\it Inverse spectral results for
the AKNS systems with partial information on the potentials}, Math.
Phys. Anal. Geom., 4 (2001), no.3, 229-244.

[GS] F. Gesztesy and B. Simon, {\it Inverse spectral analysis with
partial information on the potential. II. The case of discrete
spectrum}, Trans. Amer. Math. Soc. 352 (2000), no.6, 2765-2787.

[HH] H. Hochstadt and B. Lieberman, {\it An inverse Sturm-Liouville
problem with mixed given data}, SIAM J. Appl. Math. 34 (1978),
676-680.

[L] B. J. Levin, {\it Distribution of zeros of entire functions},
Trans. math. Mon. AMS 5 (1964).

[LG] B. Levitan and M. G. Gasymov, {\it Determination of a
differential equation by two of its spectra}, Russ. Math. Surv. 19
(1964), no.2, 1-63.

\end{document}